\documentclass[11pt]{article}
\usepackage[a4paper,margin=30mm]{geometry}
\usepackage[T1]{fontenc}

\usepackage{amsmath,amsfonts,amssymb,amsthm,mathrsfs,esint}
\usepackage{enumerate}

\usepackage{url} 
\usepackage[colorlinks=true, linkcolor=blue, citecolor=blue, urlcolor=black]{hyperref}

\numberwithin{equation}{section}
\newtheorem{theorem}{Theorem}[section]
\newtheorem{lemma}[theorem]{Lemma}

\theoremstyle{definition}
\newtheorem{definition}[theorem]{Definition}
\theoremstyle{remark}

\allowdisplaybreaks

\title
{A short route to Diening's averaging criterion\\ in spaces of homogeneous type}

\author{Alina Shalukhina
\thanks{This work is funded by national funds through 
the FCT -- Funda\c{c}\~{a}o para a Ci\^encia e a Tecnologia, I.P., 
under the scope of the projects \texttt{UID/00297/2025} and 
\texttt{UID/PRR/00297/2025} (Center for Mathematics and Applications -- NOVA Math) 
and under the scope of the PhD scholarship \texttt{UI/BD/154284/2022}.} 
\\[.5ex]
\small Center for Mathematics and Applications and Department of Mathematics, \\
\small NOVA School of Science and Technology, NOVA University Lisbon \\
\small Quinta da Torre, Caparica, 2829--516, Portugal \\
\small \texttt{a.shalukhina@campus.fct.unl.pt}
}

\begin{document}
\date{}
\maketitle

\begin{abstract}
\noindent
We present a short, self-contained proof of the averaging-operator criterion
for the boundedness of the Hardy--Littlewood maximal operator on variable Lebesgue 
spaces over unbounded spaces of homogeneous type. Building on techniques of Diening 
and Lerner, we establish a Jensen-type inequality with uniformly summable error, from 
which the criterion follows directly via a Calder\'on--Zygmund argument.

\vspace{1em}
\\[0.5em]
\textbf{MSC~2020:} 42B25, 46E30, 30L99.
\end{abstract}

\section{Introduction}

In 2005, Diening~\cite{LD05} characterized the boundedness of the Hardy--Littlewood 
maximal operator~$M$ on variable Lebesgue spaces $L^{p(\cdot)}(\mathbb{R}^n)$, for 
$1<p_-\le p_+<\infty$, by the uniform boundedness of averaging operators 
associated with disjoint cubes. The result has far-reaching consequences: it enables 
Rubio de Francia extrapolation in $L^{p(\cdot)}$ spaces~\cite{CFMP06} and 
is central to Lerner's exponent-based boundedness criteria~\cite{AL25,AL26}. 

Recently, we extended Diening's criterion to \emph{unbounded} spaces of homogeneous 
type by replacing Euclidean cubes with their abstract counterparts 
from the adjacent dyadic grids $\mathscr{D}=\bigcup_{t=1}^K\mathscr{D}^t$
of Hyt\"onen and Kairema~\cite{HK12}. Specifically, we proved the following 
theorem~\cite[Theorem~1.2]{AS26}:
\begin{theorem}[Averaging criterion]\label{th:main}
Let $(\Omega,d,\mu)$ be an \emph{unbounded} space of homogeneous type with a 
\emph{Borel-semiregular measure}~$\mu$, and let~$\mathscr{D}$ be a dyadic system 
on~$\Omega$. Suppose $p(\cdot)\in\mathcal{P}(\Omega,\mu)$ satisfies 
$1<p_-\le p_+<\infty$. Then~$M$ is bounded on~$L^{p(\cdot)}(\Omega,\mu)$ 
\emph{if and only if} the averaging operators~$T_\mathcal{Q}$, defined by
\[
T_\mathcal{Q}f:=\sum_{Q\in\mathcal{Q}} 
\bigg(\frac{1}{\mu(Q)}\int_Q|f(y)|\,d\mu(y)\bigg)\,\mathbf{1}_Q
\]
for all families~$\mathcal{Q}$ of pairwise disjoint cubes from~$\mathscr{D}$,
are uniformly bounded on~$L^{p(\cdot)}(\Omega,\mu)$.
\end{theorem}
Both Diening's original proof in~$\mathbb{R}^n$~\cite{LD05,DHHR11}
and its generalization in~\cite{AS26} are exceptionally long and technically
demanding. This paper provides a short, direct, and self-contained
proof of Theorem~\ref{th:main} that entirely bypasses the abstract Musielak--Orlicz 
machinery. The novelty of our work therefore lies in the proof method, not in 
the characterization. After collecting preliminaries in 
Section~\ref{sec:prelim}, our argument proceeds in three concise steps:
\begin{enumerate}
\item In Section~\ref{sec:RHI}, we adapt Diening's level-set norm decay 
technique~\cite[Lemma~5.5]{LD05} directly to $L^{p(\cdot)}(\Omega,\mu)$ to establish 
reverse H\"older and self-improving properties of the averaging condition.

\item In Section~\ref{sec:jensen}, we combine these properties with Lerner's dynamic 
parameter and threshold techniques~\cite{AL10,AL17} 
to develop an analytical tool of independent interest: a
\emph{Jensen-type inequality with summable error}, which extends the classical 
inequality of Kopaliani~\cite[Lemma~2.4]{TK07}.

\item In Section~\ref{sec:proof}, we prove Theorem~\ref{th:main} via 
a Calder\'on--Zygmund decomposition and the application of our Jensen-type inequality 
to rescaled, level-dependent truncations of the function.
\end{enumerate}

Note that Theorem~\ref{th:main} focuses on unbounded spaces, as the bounded
case is much simpler. For a \emph{bounded} space of homogeneous type with a 
Borel-semiregular measure, we show in~\cite[Corollary~7.4]{AS26ends} that, 
assuming $p_->1$ and permitting $p_+=\infty$, the maximal operator~$M$
is bounded on~$L^{p(\cdot)}(\Omega,\mu)$ if and only if the 
Muckenhoupt condition holds. Since dyadic cubes and balls 
are geometrically interchangeable, this takes the form
\[
\sup_{Q\in\mathscr{D}}\frac{1}{\mu(Q)}\|\mathbf{1}_Q\|_{p(\cdot)}
\|\mathbf{1}_Q\|_{p'(\cdot)}<\infty.
\]
As in~\cite[Theorem~4.5.7]{DHHR11}, the condition is equivalent to the uniform 
boundedness of the averaging operators~$T_{\{Q\}}$ over \emph{single-cube} families
on~$L^{p(\cdot)}(\Omega,\mu)$. Thus, the disjoint-family condition of 
Theorem~\ref{th:main} simplifies to single cubes 
in the bounded setting, completing the picture across both geometric regimes.

\section{Preliminaries}
\label{sec:prelim}
Variable Lebesgue spaces~$L^{p(\cdot)}(\Omega,\mu)$ are defined over 
general $\sigma$-finite, complete measure spaces~$(\Omega,\mu)$. 
Denote by $L^0(\Omega,\mu)$ the space of all $\mu$-measurable functions on~$\Omega$. 
Let $\mathcal{P}(\Omega,\mu)\subseteq L^0(\Omega,\mu)$ be the set of all 
\emph{variable exponents} $p(\cdot)\,:\,\Omega\to(1,\infty)$. Throughout, we focus on 
exponents satisfying $1<p_-\le p_+<\infty$, where
\[
p_-:=\operatornamewithlimits{ess\,inf}_{x\in\Omega}p(x)
\quad\text{and}\quad
p_+:=\operatornamewithlimits{ess\,sup}_{x\in\Omega}p(x).
\]
The conjugate exponents $p'(\cdot)\in\mathcal{P}(\Omega,\mu)$ are defined pointwise by 
$\frac{1}{p(\cdot)}+\frac{1}{p'(\cdot)}=1$.

Each~$p(\cdot)\in\mathcal{P}(\Omega,\mu)$ induces the 
\emph{modular}~$\varrho_{p(\cdot)}$, defined for all $f\in L^0(\Omega,\mu)$ by 
\[
\varrho_{p(\cdot)}(f):=\int_\Omega |f(x)|^{p(x)}\,d\mu(x).
\]
The \emph{variable Lebesgue space} $L^{p(\cdot)}(\Omega,\mu)$ then consists of 
all functions~$f\in L^0(\Omega,\mu)$ with finite Luxemburg norm
\[
\|f\|_{p(\cdot)}
:=\inf\big\{\lambda>0\,:\,\varrho_{p(\cdot)}\big(\tfrac{f}{\lambda}\big)\le1\big\}
<\infty.
\]
Equipped with this norm, $L^{p(\cdot)}(\Omega,\mu)$ is a Banach function 
space~\cite[Theorem~3.2.13]{DHHR11} that reduces to the classical Lebesgue space 
when~$p(\cdot)$ is constant. We recall basic properties of variable Lebesgue spaces.
\begin{theorem}[Toolbox for variable exponents]\label{th:toolbox}
Let $p(\cdot)\in\mathcal{P}(\Omega,\mu)$. The following hold:
\begin{enumerate}
\item[\textnormal{(a)}] \emph{(unit ball property)} $\|f\|_{p(\cdot)}\le1$
if and only if $\varrho_{p(\cdot)}(f)\le1$.

\item[\textnormal{(b)}] \emph{(unit sphere property)} If $p_+<\infty$, then 
$\|f\|_{p(\cdot)}=1$ is equivalent to $\varrho_{p(\cdot)}(f)=1$.

\item[\textnormal{(c)}] \emph{(H\"older's inequality)} For all 
$f\in L^{p(\cdot)}(\Omega,\mu)$ and $g\in L^{p'(\cdot)}(\Omega,\mu)$,
\[
\int_\Omega |f(x)\,g(x)|\,d\mu(x)\le2\,\|f\|_{p(\cdot)}\|g\|_{p'(\cdot)}.
\]

\item[\textnormal{(d)}] \emph{(norm conjugate formula)} For all $f\in L^0(\Omega,\mu)$,
\[
\|f\|_{p(\cdot)}\le2\sup_{\|g\|_{p'(\cdot)}\le1}
\int_\Omega |f(x)\,g(x)|\,d\mu(x).
\]

\item[\textnormal{(e)}] \emph{(modular-to-norm bound)} Let $p_+<\infty$. If 
$\varrho_{p(\cdot)}(f)\le ab$ for some $a\ge1$ and $b\in(0,1]$, then 
$\|f\|_{p(\cdot)}\le ab^\frac{1}{p_+}$. In particular, if $\varrho_{p(\cdot)}(f)\le1$,
then $\|f\|_{p(\cdot)}\le\varrho_{p(\cdot)}(f)^\frac1{p_+}$.

\item[\textnormal{(f)}] \emph{(convexity scaling)} For all 
$f\in L^0(\Omega,\mu)$ and $\lambda>1$, one has
$\varrho_{p(\cdot)}(\lambda f)\ge\lambda\,\varrho_{p(\cdot)}(f)$.

\item[\textnormal{(g)}] \emph{(Fatou property)} For $f_k,f\in L^0(\Omega,\mu)$,
if $0\le f_k\nearrow f$ $\mu$-almost everywhere, then
\[
\|f\|_{p(\cdot)}=\lim_{k\to\infty}\|f_k\|_{p(\cdot)}.
\]

\item[\textnormal{(h)}] \emph{(Riesz--Fischer property)} For any nonnegative sequence
$\{f_k\}_{k=0}^\infty\subseteq L^0(\Omega,\mu)$,
\[
\bigg\|\sum_{k=0}^\infty f_k\bigg\|_{p(\cdot)}\le\sum_{k=0}^\infty \|f_k\|_{p(\cdot)}.
\]
\end{enumerate}
\end{theorem}
\begin{proof}
Most of the properties can be found in the monograph~\cite{DHHR11}:
(a)~and~(b) correspond to Lemma~2.1.14; (c)~is part of Lemma~3.2.20; 
(d)~is Corollary~3.2.14; (f)~is contained in formula~(2.1.5); (g)~follows
from Theorem~2.3.17. Property~(e) is easily checked by the definition of the 
Luxemburg norm. The Riesz--Fischer property~(h) is deduced from the Fatou
property~(g) exactly as in~\cite[Chapter~1, Theorem~1.6]{BS88}.
\end{proof}

We next specialize the underlying space~$(\Omega,\mu)$ to the spaces of homogeneous 
type of Coifman and Weiss~\cite[Section~2]{CW77}. 
\begin{definition}
A \emph{space of homogeneous type} is a triple $(\Omega,d,\mu)$ with 
structural constants $\kappa\ge1$ and $C_\mu>1$ such that: 
\begin{enumerate}
\item[\textnormal{(i)}] $d$ is a \emph{quasi-metric} on~$\Omega$ satisfying 
\[
d(x,z)\le\kappa\,(d(x,y)+d(y,z))
\quad\text{for all }x,y,z\in\Omega;
\]

\item[\textnormal{(ii)}] $\mu$ is a complete Borel measure on~$\Omega$ such that 
balls $B(x,r):=\{y\in\Omega\,:\,d(x,y)<r\}$ are $\mu$-measurable and satisfy the
\emph{doubling condition}
\[
0<\mu(B(x,2r))\le C_\mu\,\mu(B(x,r))<\infty.
\]
\end{enumerate}
\end{definition}
The space $(\Omega,d,\mu)$ is \emph{bounded} if and only if
$\mu(\Omega)<\infty$~\cite[Lemma~1.9]{BC96}. Following Alvarado and 
Mitrea~\cite[Definition~2.9]{AM15}, we interpret a Borel measure as one whose 
$\sigma$-algebra of $\mu$-measurable sets \emph{contains} the Borel $\sigma$-algebra 
generated by the quasi-metric topology. This motivates the following regularity 
condition.
\begin{definition}
A space of homogeneous type~$(\Omega,d,\mu)$ has a
\emph{Borel-semiregular} measure~$\mu$ if, for every $E\subseteq\Omega$ with 
$\mu(E)<\infty$, there exists a Borel set $B\subseteq\Omega$ which differs from~$E$ 
by a set of measure zero.
\end{definition}
The local space $L^1_\mathrm{loc}(\Omega,\mu)$ consists of all 
$f\in L^0(\Omega,\mu)$ such that $f\mathbf{1}_B\in L^1(\Omega,\mu)$ for every 
ball $B\subseteq\Omega$. For any set $E\subseteq\Omega$ with $0<\mu(E)<\infty$, we 
denote the average of $f\in L^0(\Omega,\mu)$ by 
\[
|f|_E:=\fint_E |f(y)|\,d\mu(y)=\frac{1}{\mu(E)}\int_E |f(y)|\,d\mu(y).
\]
The \emph{Hardy--Littlewood maximal operator}~$M$ is defined on~$L^0(\Omega,\mu)$ by
\[
Mf(x):=\sup_{B\ni x}|f|_B, 
\quad x\in\Omega,
\]
where the supremum is taken over all balls $B\subseteq\Omega$ containing~$x$.

There is a dyadic structure on spaces of homogeneous type. The following construction of 
dyadic grids is due to Hyt\"onen and 
Kairema~\cite[Theorems~2.2 and~4.1, Corollary~7.4]{HK12}, presented here in a simplified 
form with details from~\cite[Theorem~3.4]{AS26}.
\begin{theorem}[Dyadic system]\label{th:HK}
Let $(\Omega,d,\mu)$ be a space of homogeneous type. There exist constants 
$C_\mathrm{D},C\ge1$, depending on the structural constants of the space, and 
a finite system 
\[
\mathscr{D}:=\bigcup_{t=1}^K \mathscr{D}^t
\]
of dyadic grids~$\mathscr{D}^t$, each consisting of countably many Borel sets
$Q\subseteq\Omega$ called \emph{dyadic cubes}, such that:
\begin{enumerate}
\item[\textnormal{(a)}] each grid decomposes into generations,
$\mathscr{D}^t=\bigcup_{k\in\mathbb{Z}}\mathscr{D}^t_k$, where every 
subcollection~$\mathscr{D}^t_k$ is a partition of~$\Omega$;

\item[\textnormal{(b)}] two dyadic cubes in~$\mathscr{D}^t$ are either disjoint or one 
contains the other;

\item[\textnormal{(c)}] every cube $Q\in\mathscr{D}_k^t$ has a unique \emph{parent} 
in~$\mathscr{D}_{k-1}^t$ containing it, and at least one \emph{child}
in~$\mathscr{D}_{k+1}^t$ contained in it;

\item[\textnormal{(d)}] if $Q'\in\mathscr{D}_{k+1}^t$ is a child of
$Q\in\mathscr{D}_k^t$, then 
\[
\mu(Q)\le C_\mathrm{D}\,\mu(Q'),
\]
where $C_\mathrm{D}$ is called the \emph{dyadic doubling constant};

\item[\textnormal{(e)}] every cube $Q\in\mathscr{D}$ is contained in a ball 
$B_Q\subseteq\Omega$ with $\mu(B_Q)\le C\,\mu(Q)$; conversely, every ball 
$B\subseteq\Omega$ is contained in a cube $Q_B\in\mathscr{D}$ with 
$\mu(Q_B)\le C\,\mu(B)$.
\end{enumerate}
\end{theorem}
Each dyadic grid~$\mathscr{D}^t$ induces a \emph{dyadic maximal 
operator}~$M^{\mathscr{D}^t}$ on $L^0(\Omega,\mu)$, defined~by
\[
M^{\mathscr{D}^t}f(x):=\sup_{\substack{Q\ni x,\\Q\in\mathscr{D}^t}}|f|_Q,
\quad x\in\Omega,
\]
and we set
\[
M^\mathscr{D}f(x):=\max_{1\le t\le K} M^{\mathscr{D}^t}f(x)=
\sup_{\substack{Q\ni x,\\Q\in\mathscr{D}}}|f|_Q.
\]
\begin{lemma}\label{le:M-MD}
Let $(\Omega,d,\mu)$ be a space of homogeneous type, and 
$\mathscr{D}:=\bigcup_{t=1}^K \mathscr{D}^t$ be a dyadic system on~$\Omega$.
Let $C\ge1$ be the constant from Theorem~\ref{th:HK}\textnormal{(e)}. Then: 
\begin{enumerate}
\item[\textnormal{(a)}] \emph{(equivalence of operators)} 
For all $f\in L^0(\Omega,\mu)$,
\[
\frac{1}{C}\,Mf(x)\le M^\mathscr{D}f(x)\le C\,Mf(x)
\quad\text{for every }x\in\Omega.
\]

\item[\textnormal{(b)}] \emph{(majorant)} If $\mu$ is Borel-semiregular, then for all 
$t\in\{1,\ldots,K\}$ and $f\in L^1_\mathrm{loc}(\Omega,\mu)$,
\[
|f(x)|\le M^{\mathscr{D}^t}f(x)
\quad\text{for $\mu$-almost every }x\in\Omega.
\]
\end{enumerate}
\end{lemma}
\begin{proof}
Equivalence~(a) is a direct consequence of Theorem~\ref{th:HK}(e). 
Part~(b) is established in~\cite[Lemma~3.10]{AS26}.
\end{proof}

The level sets of each dyadic maximal function $M^{\mathscr{D}^t}f$ admit the following 
Calder\'on--Zygmund decomposition~\cite[Lemma~3.12]{AS26}.
\begin{lemma}[Calder\'on--Zygmund decomposition]\label{le:CZ}
Let $(\Omega,d,\mu)$ be an \emph{unbounded} space of homogeneous type, and let 
$\mathscr{D}=\bigcup_{t=1}^K \mathscr{D}^t$ be a dyadic system on~$\Omega$ 
with dyadic doubling constant $C_\mathrm{D}\ge1$.
Given $f\in L^1(\Omega,\mu)$, $\lambda>0$, and $t\in\{1,\ldots,K\}$, there 
exists a family $\mathcal{Q}$
of \emph{pairwise disjoint} cubes $Q\in\mathscr{D}^t$ with
$\lambda<|f|_Q\le C_\mathrm{D}\,\lambda$ such that
\[
\Omega_\lambda^t:=\left\{x\in\Omega\,:\,
M^{\mathscr{D}^t}f(x)>\lambda
\right\}=\bigcup_{Q\in\mathcal{Q}}Q.
\]
Moreover, the cubes $Q\in\mathcal{Q}$ are \emph{maximal} in~$\Omega_\lambda^t$: no cube 
from~$\mathscr{D}^t$ strictly containing some $Q\in\mathcal{Q}$ is contained 
in~$\Omega_\lambda^t$.
\end{lemma}
Finally, let~$\mathfrak{S}$ denote the set of all families~$\mathcal{Q}$ of pairwise
disjoint cubes from a dyadic system~$\mathscr{D}$ on~$\Omega$. For each 
$\mathcal{Q}\in\mathfrak{S}$, the \emph{averaging operator} 
$T_\mathcal{Q}\,:\,L^1_\mathrm{loc}(\Omega,\mu)\to L^0(\Omega,\mu)$
is defined by 
\[
T_\mathcal{Q}f(x):=\sum_{Q\in\mathcal{Q}} |f|_Q\,\mathbf{1}_Q(x),
\quad x\in\Omega.
\]
The following condition is necessary for the boundedness of the maximal 
operator~$M$ on~$L^{p(\cdot)}(\Omega,\mu)$.
\begin{definition}
Given a dyadic system~$\mathscr{D}$ on~$\Omega$, an exponent 
$p(\cdot)\in\mathcal{P}(\Omega,\mu)$ satisfies the 
\emph{strong dyadic Muckenhoupt condition}, denoted 
$p(\cdot)\in\mathcal{A}^\mathscr{D}(\Omega,\mu)$, if
\[
[p]_\mathcal{A}^\mathscr{D}:=\sup_{\mathcal{Q}\in\mathfrak{S}}
\|T_\mathcal{Q}\|_{L^{p(\cdot)}\to L^{p(\cdot)}}<\infty.
\]
\end{definition}
Indeed, necessity follows from the pointwise bound
\begin{equation}\label{ineq:necessity}
T_\mathcal{Q}f(x)\le M^\mathscr{D}f(x)\le C\,Mf(x),
\quad x\in\Omega,
\end{equation}
for all $\mathcal{Q}\in\mathfrak{S}$ and $f\in L^1_\mathrm{loc}(\Omega,\mu)$,
where the first inequality holds by definition and the second is due to
Lemma~\ref{le:M-MD}(a). Conversely, the rest of the paper establishes that 
$p(\cdot)\in\mathcal{A}^\mathscr{D}(\Omega,\mu)$ is also sufficient 
when $1<p_-\le p_+<\infty$, provided the space $(\Omega,d,\mu)$ is unbounded 
and its measure~$\mu$ is Borel-semiregular.

Note that $[p]_\mathcal{A}^\mathscr{D}\ge1$, since 
$T_{\{Q\}}(\mathbf{1}_Q)=\mathbf{1}_Q$ for all $Q\in\mathscr{D}$.
We conclude this section by recording the duality of the 
class~$\mathcal{A}^\mathscr{D}(\Omega,\mu)$.
\begin{lemma}[Duality of~$\mathcal{A}^\mathscr{D}$]\label{le:duality}
If $p(\cdot)\in\mathcal{A}^\mathscr{D}(\Omega,\mu)$, then 
$p'(\cdot)\in\mathcal{A}^\mathscr{D}(\Omega,\mu)$ with 
$[p']_\mathcal{A}^\mathscr{D}\le 4[p]_\mathcal{A}^\mathscr{D}$.
\end{lemma}
\begin{proof}
For every~$\mathcal{Q}\in\mathfrak{S}$ and all~$f,g\in L^1_\mathrm{loc}(\Omega,\mu)$, 
we have the symmetry relation
\[
\int_\Omega T_\mathcal{Q}f(x)\,|g(x)|\,d\mu(x)
=\sum_{Q\in\mathcal{Q}}\mu(Q)\,|f|_Q\,|g|_Q
=\int_\Omega T_\mathcal{Q}g(x)\,|f(x)|\,d\mu(x).
\]
The norm conjugate formula and H\"older's inequality then yield
\[
\|T_\mathcal{Q}\|_{L^{p'(\cdot)}\to L^{p'(\cdot)}}
\le 2\sup_{\|f\|_{p(\cdot)}\le1}\sup_{\|g\|_{p'(\cdot)}\le1}
\int_\Omega T_\mathcal{Q}f(x)\,|g(x)|\,d\mu(x)
\le 4\,\|T_\mathcal{Q}\|_{L^{p(\cdot)}\to L^{p(\cdot)}},
\]
and taking the supremum over $\mathcal{Q}\in\mathfrak{S}$ gives the claim.
\end{proof}

\section{Reverse H\"older Property and Self-Improvement of $\mathcal{A}^\mathscr{D}$}
\label{sec:RHI}

From now on, we fix an \emph{unbounded} space of homogeneous 
type~$(\Omega,d,\mu)$ with a \emph{Borel-semiregular measure}~$\mu$, and 
a dyadic system~$\mathscr{D}:=\bigcup_{t=1}^K \mathscr{D}^t$ from Theorem~\ref{th:HK}.

We first establish reverse H\"older and self-improving properties for
$p(\cdot)\in\mathcal{A}^\mathscr{D}(\Omega,\mu)$ using Diening's 
level-set norm decay technique from~\cite[Lemma~5.5]{LD05}. Claim~(a) recovers
the reverse H\"older results by Lerner~\cite[Lemma~5.1]{AL17} and 
Karlovych~\cite[Lemma~5.1]{AK20}, but without their \emph{a priori} assumption that~$M$ 
is bounded on~$L^{p(\cdot)}$. Claim~(b) is a variable-Lebesgue analogue 
of~\cite[Theorem~5.6]{AS26}.
\begin{lemma}\label{le:RHI}
Let $p(\cdot)\in\mathcal{A}^\mathscr{D}(\Omega,\mu)$ with $1<p_-\le p_+<\infty$. 
There exist a parameter $r\in(1,p_-)$ and constants $R,C_r\ge1$ such that 
the following hold:
\begin{enumerate}
\item[\textnormal{(a)}] \emph{(reverse H\"older property)}
For every family $\mathcal{Q}\in\mathfrak{S}$ and every nonnegative 
sequence $\{t_Q\}_{Q\in\mathcal{Q}}$,
\[
\sum_{Q\in\mathcal{Q}} \int_Q t_Q^{p(x)}\,d\mu(x)\le1
\ \implies\ 
\sum_{Q\in\mathcal{Q}} \mu(Q) \bigg(\fint_Q t_Q^{rp(x)}\,d\mu(x)\bigg)^\frac1r\le R.
\]
The same implication holds for the conjugate exponent $p'(\cdot)$.

\item[\textnormal{(b)}] \emph{(self-improvement of~$\mathcal{A}^\mathscr{D}$)}
For every family $\mathcal{Q}\in\mathfrak{S}$ and every nonnegative function 
$f \in L^{p(\cdot)}(\Omega,\mu)$,
\[
\bigg\|\sum_{Q\in\mathcal{Q}} |f^r|_Q^{\frac{1}{r}}\,\mathbf{1}_Q\bigg\|_{p(\cdot)} 
\le C_r\,\|f\|_{p(\cdot)}.
\]
Equivalently, the scaled exponent $\frac{p(\cdot)}{r}$ belongs to 
$\mathcal{A}^\mathscr{D}(\Omega,\mu)$. 
\end{enumerate}
\end{lemma}
\begin{proof}
\textbf{Step~1: Norm decay on level sets.} Take a family $\mathcal{Q}\in\mathfrak{S}$ 
and a nonnegative function $h\in L^1_{\mathrm{loc}}(\Omega,\mu)$ with $|h|_Q>0$ 
for all $Q\in\mathcal{Q}$. For a fixed cube~$Q\in\mathcal{Q}$, find 
$t\in\{1,\ldots,K\}$ such that $Q\in\mathscr{D}^t$. 
We consider the sets $\Omega_0^Q:=Q$ and, for $k\ge1$,
\[
\Omega_k^Q:=\left\{x\in\Omega\,:\,
M^{\mathscr{D}^t}(h\mathbf{1}_Q)(x)>a^k\,|h|_Q
\right\}
\]
with $a:=2C_\mathrm{D}$, where~$C_\mathrm{D}$ is the dyadic doubling constant 
of~$\mathscr{D}$. Note that $\Omega_k^Q\subseteq Q$ for all $k\ge0$.
Indeed, if $x\notin Q$, any dyadic cube $P\in\mathscr{D}^t$ containing~$x$ with 
$|h\mathbf{1}_Q|_P>0$ must strictly contain~$Q$. Since $\mu(Q)<\mu(P)$, we have
$|h\mathbf{1}_Q|_P<|h|_Q$, and hence  
$M^{\mathscr{D}^t}(h\mathbf{1}_Q)(x)\le |h|_Q<a^k\,|h|_Q$ 
for all $x\notin Q$ and $k\ge1$, ensuring $\Omega_k^Q\subseteq Q$. 

By the Calder\'on--Zygmund Lemma~\ref{le:CZ} applied 
to~$h\mathbf{1}_Q\in L^1(\Omega,\mu)$, each nonempty~$\Omega_k^Q$ for $k\ge1$ 
decomposes into a family 
$\mathcal{P}_k^Q\in\mathfrak{S}$ of pairwise disjoint maximal cubes 
$V\in\mathscr{D}^t$ such that
\begin{equation}\label{est:threshold}
a^k\,|h|_Q<|h|_V\le C_\mathrm{D}\,a^k\,|h|_Q.
\end{equation}
For $k=0$, set $\mathcal{P}_0^Q:=\{Q\}$, so the upper bound in~\eqref{est:threshold} 
holds trivially since $C_\mathrm{D}\ge1$. Every cube $W\in\mathcal{P}_{k+1}^Q$ 
intersecting $V\in\mathcal{P}_k^Q$ is contained in~$V$; otherwise, 
$V\subsetneq W\subseteq\Omega_{k+1}^Q\subseteq\Omega_k^Q$ would contradict the
maximality of~$V$. For all
$V\in\mathcal{P}_k^Q$ with $k\ge0$, estimates~\eqref{est:threshold} thus yield
\[
\mu(V\cap\Omega_{k+1}^Q)
=\sum_{\substack{W\in\mathcal{P}_{k+1}^Q:\\W\subseteq V}}\mu(W)
\le\frac{1}{a^{k+1}\,|h|_Q}\int_V h(x)\,d\mu(x)
\le\frac{C_\mathrm{D}\,a^k\,|h|_Q}{a^{k+1}\,|h|_Q}\,\mu(V)
=\frac{1}{2}\,\mu(V),
\]
whence $\mu(V\setminus\Omega_{k+1}^Q)\ge\frac12\,\mu(V)$. Since the cubes 
in~$\mathcal{Q}$ are pairwise disjoint, the combined collection 
$\mathcal{P}_k:=\bigcup_{Q\in\mathcal{Q}}\mathcal{P}_k^Q$ belongs to~$\mathfrak{S}$ 
for each $k\ge0$. Consequently, for any nonnegative sequence 
$\{t_Q\}_{Q\in\mathcal{Q}}$, we have the pointwise estimate
\[
T_{\mathcal{P}_k}\bigg(\sum_{Q\in\mathcal{Q}}
t_Q\mathbf{1}_{\Omega_k^Q\setminus\Omega_{k+1}^Q}\bigg)
=\sum_{Q\in\mathcal{Q}} \sum_{V\in\mathcal{P}_k^Q}
\frac{\mu(V\setminus\Omega_{k+1}^Q)}{\mu(V)}\,t_Q\mathbf{1}_V
\ge\frac12 \sum_{Q\in\mathcal{Q}} t_Q\mathbf{1}_{\Omega_k^Q}.
\]

Set $\lambda_k:=\|\sum_{Q\in\mathcal{Q}}t_Q\mathbf{1}_{\Omega_k^Q}\|_{p(\cdot)}$. 
Whenever $\lambda_k>0$, the condition
$p(\cdot)\in\mathcal{A}^\mathscr{D}(\Omega,\mu)$ implies
\[
\bigg\|\sum_{Q\in\mathcal{Q}}\frac{t_Q}{\lambda_k}\,
\mathbf{1}_{\Omega_k^Q\setminus\Omega_{k+1}^Q}\bigg\|_{p(\cdot)}
\ge\frac{1}{2[p]_\mathcal{A}^\mathscr{D}},
\]
so the corresponding modular is bounded from 
below by~$(2[p]_\mathcal{A}^\mathscr{D})^{-p_+}$ by the special case of 
Theorem~\ref{th:toolbox}(e). Together with the unit sphere property, this yields
\[
\varrho_{p(\cdot)}\bigg(\sum_{Q\in\mathcal{Q}} \frac{t_Q}{\lambda_k}\, 
\mathbf{1}_{\Omega_{k+1}^Q} \bigg)
=1-\varrho_{p(\cdot)}\bigg(\sum_{Q\in\mathcal{Q}}\frac{t_Q}{\lambda_k}\, 
\mathbf{1}_{\Omega_k^Q\setminus\Omega_{k+1}^Q} \bigg) 
\le 1 - \frac1{(2[p]_\mathcal{A}^\mathscr{D})^{p_+}},
\]
whence Theorem~\ref{th:toolbox}(e) again implies 
$\lambda_{k+1}\le\beta_p\,\lambda_k$ with 
$\beta_p:=(1-(2[p]_\mathcal{A}^\mathscr{D})^{-p_+})^{\frac1{p_+}}<1$. 
By induction, for all $k\ge0$,
\begin{equation}\label{ineq:norm-decay}
\bigg\|\sum_{Q\in\mathcal{Q}} t_Q\mathbf{1}_{\Omega_k^Q}\bigg\|_{p(\cdot)}
\le\beta_p^k\,\bigg\|\sum_{Q\in\mathcal{Q}} t_Q\mathbf{1}_Q\bigg\|_{p(\cdot)}.
\end{equation}
Since $p'(\cdot)\in\mathcal{A}^\mathscr{D}(\Omega,\mu)$ by Lemma~\ref{le:duality}, 
the exact same bound holds for~$p'(\cdot)$ with constant $\beta_{p'}<1$. Set 
$\beta:=\max\{\beta_p,\beta_{p'}\}$, choose $r\in(1,\min\{p_-,(p')_-\})$ 
sufficiently close to~$1$ so that $a^{r-1}\beta<1$, and define
\[
R:=\frac{a^{r-1}}{1-a^{r-1}\beta}.
\]

Note two useful estimates. First, since $|h^r|_Q\ge|h|_Q^r$ by Jensen's inequality, 
whenever $h\in L^r(Q,\mu)$ we have
\begin{equation}\label{est:jensen-r}
|h^r|_Q^{\frac{1}{r}}
=|h^r|_Q^{\frac{1}{r}-1}\,|h^r|_Q
\le\fint_Q h(x)\bigg(\frac{h(x)}{|h|_Q}\bigg)^{r-1}\,d\mu(x).
\end{equation}
Second, for all $k\ge0$ and $\mu$-almost every $x\in\Omega_k^Q\setminus\Omega_{k+1}^Q$,
the ratio in the parentheses is bounded in view of Lemma~\ref{le:M-MD}(b) by
\begin{equation}\label{est:ratio}
\frac{h(x)}{|h|_Q}\le\frac{M^{\mathscr{D}^t}(h\mathbf{1}_Q)(x)}{|h|_Q}
\le a^{k+1}.
\end{equation}

\paragraph*{Step~2: Proof of~(a).} Fix a family $\mathcal{Q}\in\mathfrak{S}$ and a 
nonnegative sequence $\{t_Q\}_{Q\in\mathcal{Q}}$ satisfying
\[
\varrho_{p(\cdot)}\bigg(\sum_{Q\in\mathcal{Q}}t_Q\mathbf{1}_Q\bigg)
=\sum_{Q\in\mathcal{Q}}\int_Q t_Q^{p(x)}\,d\mu(x)\le1.
\]
Without loss of generality, we may assume $t_Q>0$ for all $Q\in\mathcal{Q}$.
Apply \emph{Step~1} to the function 
$h(x):=\sum_{Q\in\mathcal{Q}} t_Q^{p(x)}\,\mathbf{1}_Q(x)$, $x\in\Omega$. Note that 
$\big\|\sum_{Q\in\mathcal{Q}}t_Q\mathbf{1}_Q\big\|_{p(\cdot)}\le1$. Then the norm-decay 
estimate~\eqref{ineq:norm-decay}, the unit ball property, and the convexity scaling 
from Theorem~\ref{th:toolbox}(f) yield
\[
\sum_{Q\in\mathcal{Q}}\int_{\Omega_k^Q} h(x)\,d\mu(x)
=\varrho_{p(\cdot)}\bigg(\sum_{Q\in\mathcal{Q}}t_Q\mathbf{1}_{\Omega_k^Q}\bigg)
\le\beta^k
\]
for all $k\ge0$. Applying inequality~\eqref{est:jensen-r} together with 
the bound~\eqref{est:ratio} on the ``annuli'' $\Omega_k^Q\setminus\Omega_{k+1}^Q$
splitting~$Q$ gives
\begin{align*}
\sum_{Q\in\mathcal{Q}} \mu(Q) \bigg(\fint_Q t_Q^{rp(x)}\,d\mu(x)\bigg)^\frac1r
&\le\sum_{Q\in\mathcal{Q}} \int_Q h(x)\bigg(\frac{h(x)}{|h|_Q}\bigg)^{r-1}\,d\mu(x)
\\
&\le a^{r-1} \sum_{k=0}^\infty a^{(r-1)k}
\sum_{Q\in\mathcal{Q}} \int_{\Omega_k^Q} h(x)\,d\mu(x)
\\
&\le a^{r-1} \sum_{k=0}^\infty \big(a^{r-1}\beta\big)^k=R.
\end{align*}
The identical argument holds for~$p'(\cdot)$ with the same~$r$ and~$R$. Thus~(a) is 
proved.  

\paragraph*{Step~3: Proof of~(b).} Fix a family $\mathcal{Q}\in\mathfrak{S}$ and a 
function $f\ge0$ from $L^{p(\cdot)}(\Omega,\mu)$ with $|f|_Q>0$ for all 
$Q\in\mathcal{Q}$. Apply \emph{Step~1} to $h:=f$. Given a nonnegative sequence 
$\{t_Q\}_{Q\in\mathcal{Q}}$, it follows from the annular bound~\eqref{est:ratio}, 
the Riesz--Fischer property, and the norm-decay estimate~\eqref{ineq:norm-decay} 
for~$p'(\cdot)$ that 
\begin{align*}
\bigg\|\sum_{Q\in\mathcal{Q}} t_Q
\bigg(\frac{f}{|f|_Q}\bigg)^{r-1}\mathbf{1}_Q\bigg\|_{p'(\cdot)}
&=\bigg\| \sum_{k=0}^\infty \sum_{Q\in\mathcal{Q}} t_Q
\bigg(\frac{f}{|f|_Q}\bigg)^{r-1}
\mathbf{1}_{\Omega_k^Q\setminus\Omega_{k+1}^Q}\bigg\|_{p'(\cdot)}
\\
&\le a^{r-1}\sum_{k=0}^\infty a^{(r-1)k}\,
\bigg\|\sum_{Q\in\mathcal{Q}} t_Q\,\mathbf{1}_{\Omega_k^Q}\bigg\|_{p'(\cdot)}
\\
&\le R\,\bigg\|\sum_{Q\in\mathcal{Q}} t_Q\,\mathbf{1}_Q\bigg\|_{p'(\cdot)}.
\end{align*}
For any nonnegative $g\in L^{p'(\cdot)}(\Omega,\mu)$ with $\|g\|_{p'(\cdot)}\le1$, 
successively applying~\eqref{est:jensen-r}, H\"older's inequality, and the above 
estimate for $t_Q:=|g|_Q$ yields
\begin{align*}
\int_\Omega \bigg[\sum_{Q\in\mathcal{Q}} 
|f^r|_Q^{\frac{1}{r}}\,\mathbf{1}_Q(x)\bigg]g(x)\,d\mu(x)
&=\sum_{Q\in\mathcal{Q}}\mu(Q)\,|f^r|_Q^{\frac{1}{r}}\,|g|_Q
\\
&\le\int_\Omega f(x)\bigg[\sum_{Q\in\mathcal{Q}}|g|_Q
\bigg(\frac{f(x)}{|f|_Q}\bigg)^{r-1}\mathbf{1}_Q(x)\bigg]d\mu(x)
\\
&\le 2R\,\|f\|_{p(\cdot)}\,
\bigg\|\sum_{Q\in\mathcal{Q}}|g|_Q\,\mathbf{1}_Q\bigg\|_{p'(\cdot)}
\\
&\le 2R\,[p']_\mathcal{A}^\mathscr{D}\,\|f\|_{p(\cdot)}.
\end{align*}
The norm conjugate formula gives the desired inequality with 
$C_r:=4R\,[p']_\mathcal{A}^\mathscr{D}$. Finally, replacing $f$ by $f^\frac1r$ and 
using the rescaling $\|h^r\|_{\frac{p(\cdot)}{r}}=\|h\|_{p(\cdot)}^r$ shows that
$\frac{p(\cdot)}{r}\in\mathcal{A}^\mathscr{D}(\Omega,\mu)$.
\end{proof}

\section{Jensen-Type Inequality with Summable Error}
\label{sec:jensen}

Using the reverse H\"older and self-improving properties established above,
we extend the Euclidean Jensen-type inequalities of Kopaliani~\cite[Lemma~2.4]{TK07} 
and Lerner~\cite[Lemma~3.1]{AL10} to spaces of homogeneous type for 
$p(\cdot)\in\mathcal{A}^\mathscr{D}(\Omega,\mu)$. Crucially, we 
avoid assuming a lower bound on~$|g|_Q$. Our inequality holds for arbitrarily small 
averages at the cost of an additional \emph{error term} with uniform summability 
properties.
\begin{theorem}[Jensen-type inequality]\label{th:jensen}
Let $p(\cdot)\in\mathcal{A}^\mathscr{D}(\Omega,\mu)$ with $1<p_-\le p_+<\infty$.
There exist constants $C>1$ and $\delta>0$ and a function $b:\mathscr{D}\to[0,\infty)$ 
satisfying 
\[
C_b:=\sup_{\mathcal{Q}\in\mathfrak{S}}
\sum_{Q\in\mathcal{Q}}b(Q)<\infty
\]
such that, for all~$Q\in\mathscr{D}$ and all~$g\in L^{p(\cdot)}(\Omega,\mu)$
with $\|g\|_{p(\cdot)}\le1$,
\begin{equation}\label{ineq:Jensen}
\int_Q (|g|_Q)^{p(x)}\,d\mu(x) 
\le C\int_Q |g(x)|^{p(x)}\,d\mu(x)
+b(Q)\,(|g|_Q)^\delta\,\mathbf{1}_{(0,1)}(|g|_Q).
\end{equation}
\end{theorem}
The proof of this inequality adapts the dynamic parameter method due to 
Lerner~\cite[Lemma~3.1]{AL10}. We prepare for the main argument with two auxiliary 
lemmas. Define 
\[
I_Q(t):=\fint_Q t^{p'(y)-1}\,d\mu(y)
\quad\text{for }t>0.
\]
\begin{lemma}[Characteristic window]\label{le:window}
Let $p(\cdot)\in\mathcal{A}^\mathscr{D}(\Omega,\mu)$ with $1<p_-\le p_+<\infty$, 
and let $r>1$ be as in Lemma~\ref{le:RHI}.
There exists $C_1>1$ such that,
for every cube $Q\in\mathscr{D}$ and every scalar
\[
t\in\bigg[\min\bigg\{1,\frac{1}{\|\mathbf{1}_Q\|_{p'(\cdot)}^r}\bigg\},
\max\bigg\{1,\frac{1}{\|\mathbf{1}_Q\|_{p'(\cdot)}^r}\bigg\}\bigg],
\]
there holds
\[
\int_Q I_Q(t)^{p(x)}\,d\mu(x)\le C_1\int_Q t^{p'(y)}\,d\mu(y).
\]
\end{lemma}
\begin{proof}
Factorize the left-hand side of the desired inequality as
\begin{align*}
\int_Q I_Q(t)^{p(x)}\,d\mu(x)
&=\int_Q \bigg(\fint_Q t^{p'(y)-1}\,d\mu(y)\bigg)^{p(x)}\,d\mu(x)
\\
&=\int_Q \bigg(\fint_Q t^{p'(y)-p'(x)}\,d\mu(y)\bigg)^{p(x)}\,t^{p'(x)}\,d\mu(x)
\\
&=\fint_Q \bigg(\fint_Q t^{p'(y)-p'(x)}\,d\mu(y)\bigg)^{p(x)-1}\,d\mu(x)
\,\int_Q t^{p'(y)}\,d\mu(y).
\end{align*}
It suffices to bound the first factor. By the unit sphere property,
$\int_Q \|\mathbf{1}_Q\|_{p(\cdot)}^{-p(x)}\,d\mu(x)=1$.
Applying the reverse H\"older property from Lemma~\ref{le:RHI}(a) to the single-cube 
family $\mathcal{Q}=\{Q\}$ with $t_Q=\|\mathbf{1}_Q\|_{p(\cdot)}^{-1}$---and analogously 
for~$p'(\cdot)$---gives
\[
\fint_Q \|\mathbf{1}_Q\|_{p(\cdot)}^{-rp(x)}\,d\mu(x)
\le\bigg(\frac{R}{\mu(Q)}\bigg)^r
\quad\text{and}\quad
\fint_Q \|\mathbf{1}_Q\|_{p'(\cdot)}^{-rp'(x)}\,d\mu(x)
\le\bigg(\frac{R}{\mu(Q)}\bigg)^r.
\]
For $t$ in the specified interval and $x,y\in Q$, the monotonicity of 
$\tau\mapsto\tau^{p'(y)-p'(x)}$ between the endpoints $1$ and 
$\|\mathbf{1}_Q\|_{p'(\cdot)}^{-r}$, combined with the reverse H\"older bound 
for~$p'(\cdot)$, yields
\begin{align*}
\fint_Q t^{p'(y)-p'(x)}\,d\mu(y)
&\le\fint_Q\bigg[1+
\bigg(\frac{1}{\|\mathbf{1}_Q\|_{p'(\cdot)}}\bigg)^{r(p'(y)-p'(x))}\bigg]\,d\mu(y)
\\
&\le1+\bigg(\frac{R\,\|\mathbf{1}_Q\|_{p'(\cdot)}^{p'(x)}}{\mu(Q)}\bigg)^r.
\end{align*}
Moreover, the condition $p(\cdot)\in\mathcal{A}^\mathscr{D}(\Omega,\mu)$ 
implies via the norm conjugate formula that
\[
\frac{1}{\mu(Q)}\|\mathbf{1}_Q\|_{p(\cdot)}\|\mathbf{1}_Q\|_{p'(\cdot)}
\le 2\,\sup_{\|f\|_{p(\cdot)}\le1} \|\,|f|_Q\,\mathbf{1}_Q\|_{p(\cdot)}
=2\,\|T_{\{Q\}}\|_{L^{p(\cdot)}\to L^{p(\cdot)}}\le 2[p]_\mathcal{A}^\mathscr{D}.
\]
Using the elementary inequality $(1+u)^{p(x)-1}\le 2^{p_+}(1+u^{p(x)-1})$, we combine 
the two preceding estimates and the reverse H\"older bound for~$p(\cdot)$ to obtain
\begin{align*}
\fint_Q \bigg(\fint_Q t^{p'(y)-p'(x)}\,&d\mu(y)\bigg)^{p(x)-1}\,d\mu(x)
\\
&\le 2^{p_+} \bigg[1+\fint_Q 
\bigg(\frac{R\,\|\mathbf{1}_Q\|_{p'(\cdot)}^{p'(x)}}
{\mu(Q)}\bigg)^{r(p(x)-1)}\,d\mu(x)\bigg]
\\
&\le 2^{p_+} \bigg[1+R^{r(p_+-1)}\,\mu(Q)^r\fint_Q 
\bigg(\frac{\|\mathbf{1}_Q\|_{p'(\cdot)}}{\mu(Q)}\bigg)^{rp(x)}\,d\mu(x)\bigg]
\\
&\le 2^{p_+} \bigg[1+(2R\,[p]_\mathcal{A}^\mathscr{D})^{rp_+}
\,\bigg(\frac{\mu(Q)}{R}\bigg)^r\fint_Q 
\|\mathbf{1}_Q\|_{p(\cdot)}^{-rp(x)}\,d\mu(x)\bigg]
\\
&\le 2^{p_+} \big[1+(2R\,[p]_\mathcal{A}^\mathscr{D})^{rp_+}\big]=:C_1,
\end{align*}
which completes the proof.
\end{proof}

Our second auxiliary lemma relates two quantities based on
the parameter~$r>1$ from Lemma~\ref{le:RHI}. For $Q\in\mathscr{D}$ and $t>0$, define
\[
U_Q(t):=\mu(Q) \bigg(\fint_Q t^{p'(x)}\,d\mu(x)\bigg)^\frac1r
\quad\text{and}\quad
V_Q(t):=\mu(Q) \bigg(\fint_Q I_Q(t)^{p(x)}\,d\mu(x)\bigg)^\frac1r.
\]
If $1<p_-\le p_+<\infty$, the functions $t\mapsto U_Q(t)$ and $t\mapsto V_Q(t)$ are 
continuous and strictly increasing, and they extend by continuity to
$U_Q(0)=V_Q(0)=0$. The following domination relation holds between their sums over 
disjoint cubes.
\begin{lemma}[Multi-cube domination]\label{le:domination}
Let $p(\cdot)\in\mathcal{A}^\mathscr{D}(\Omega,\mu)$ with $1<p_-\le p_+<\infty$,
and let $r>1$ and $R\ge1$ be as in Lemma~\ref{le:RHI}. 
There exists $C_2>1$ such that, for every family $\mathcal{Q}\in\mathfrak{S}$ and 
every nonnegative sequence $\{t_Q\}_{Q\in\mathcal{Q}}$,
\[
\sum_{Q\in\mathcal{Q}} U_Q(t_Q)\le R
\ \implies\ 
\sum_{Q\in\mathcal{Q}} 
V_Q(t_Q)\le C_2.
\]
\end{lemma}
\begin{proof}
Consider the function
\[
f(x):=\sum_{Q\in\mathcal{Q}} t_Q^\frac{p'(x)-1}{r}\,\mathbf{1}_Q(x),
\quad x\in\Omega.
\]
By Jensen's inequality, since the cubes in~$\mathcal{Q}$ are pairwise
disjoint, 
\[
\varrho_{p(\cdot)}(f)=\sum_{Q\in\mathcal{Q}} 
\int_Q t_Q^{\frac{p'(x)}{r}}\,d\mu(x)
\le\sum_{Q\in\mathcal{Q}} \mu(Q)
\bigg(\fint_Q t_Q^{p'(x)}\,d\mu(x)\bigg)^{\frac1r}
\le R.
\]
By Theorem~\ref{th:toolbox}(e), this implies $\|f\|_{p(\cdot)}\le R$. 
Moreover, $|f^r|_Q=I_Q(t_Q)$ on each~$Q\in\mathcal{Q}$, so the self-improvement bound 
from Lemma~\ref{le:RHI}(b) gives
\[
\bigg\|\sum_{Q\in\mathcal{Q}} I_Q(t_Q)^{\frac{1}{r}}\,\mathbf{1}_Q\bigg\|_{p(\cdot)}
\le C_r\,\|f\|_{p(\cdot)} \le R C_r.
\]
By the unit ball property,
\[
\sum_{Q\in\mathcal{Q}}\int_Q 
\bigg(\frac{I_Q(t_Q)^\frac{1}{r}}{RC_r}\bigg)^{p(x)}\,d\mu(x)
=\varrho_{p(\cdot)}\bigg(\sum_{Q\in\mathcal{Q}}
\frac{I_Q(t_Q)^{\frac{1}{r}}}{RC_r}\,\mathbf{1}_Q\bigg)\le1.
\]
Applying the reverse H\"older property from Lemma~\ref{le:RHI}(a) with~$t_Q$ replaced by 
$\frac{I_Q(t_Q)^{1/r}}{RC_r}$ and using $R C_r\ge1$ therefore yields
\[
\sum_{Q\in\mathcal{Q}}\mu(Q)\bigg(\fint_Q I_Q(t_Q)^{p(x)}\,d\mu(x)\bigg)^{\frac{1}{r}}
\le (RC_r)^{p_+} R =: C_2,
\]
which completes the proof.
\end{proof}

\begin{proof}[Proof of Theorem~\ref{th:jensen}]
Let $r>1$ and $R\ge1$ be as in Lemma~\ref{le:RHI},
and let $C_1,C_2>1$ be the constants from Lemmas~\ref{le:window} 
and~\ref{le:domination}. Set $\lambda:=2^{\frac{r}{(p')_--1}}>1$ and
$\gamma:=\frac{((p')_+-1)\,p_+}{r}>0$, so that for all~$Q\in\mathscr{D}$ and~$t>0$,
\[
U_Q\Big(\frac{t}{\lambda}\Big)\le\frac12\,U_Q(t)
\quad\text{and}\quad
V_Q(t)\le\lambda^\gamma\,V_Q\Big(\frac{t}{\lambda}\Big)
\quad\text{and}\quad
2\,I_Q(t)\le I_Q(\lambda t).
\]

We first construct the summable error $b(Q)$ following the threshold argument
due to Lerner~\cite[Lemma~5.2]{AL17}. Fix
\[
\kappa:=\lambda^\gamma\max\bigg\{C_1^{\frac{1}{r}},\,\frac{C_2}{R}\bigg\}.
\]
For each $Q\in\mathscr{D}$, let $t_Q\ge0$ be the supremum of the set 
\[
\mathcal{E}_Q:=\bigg\{t>0\ :\ U_Q(t)\le R 
\quad\text{and}\quad 
V_Q(t)>\kappa\, U_Q(t) \bigg\},
\]
with $t_Q:=0$ if $\mathcal{E}_Q=\emptyset$. Note that $U_Q(t_Q)<R$: otherwise, 
$U_Q(t_Q)=R$ would imply $V_Q(t_Q)\ge\kappa R>C_2$, contradicting 
Lemma~\ref{le:domination} for $\mathcal{Q}=\{Q\}$.
By continuity and the extremality of~$t_Q$, the strict inequality $U_Q(t_Q)<R$ forces 
$V_Q(t_Q)=\kappa\,U_Q(t_Q)$. Define
\[
b(Q):=V_Q(t_Q).
\]
To check~$C_b<\infty$, observe that for any family 
$\mathcal{Q}\in\mathfrak{S}$ satisfying $\sum_{Q\in\mathcal{Q}} U_Q(t_Q)\le 2R$, 
we have $\sum_{Q\in\mathcal{Q}} U_Q(\frac{t_Q}{\lambda})\le R$. 
Lemma~\ref{le:domination} then implies
\begin{equation}\label{est:b-one-family}
\sum_{Q\in\mathcal{Q}} b(Q)
=\sum_{Q\in\mathcal{Q}} V_Q(t_Q)
\le\lambda^\gamma \sum_{Q\in\mathcal{Q}} V_Q\Big(\frac{t_Q}{\lambda}\Big)
\le \lambda^\gamma C_2,
\end{equation}
and consequently
\[
\sum_{Q\in\mathcal{Q}} U_Q(t_Q)
=\frac{1}{\kappa}\sum_{Q\in\mathcal{Q}} V_Q(t_Q)
\le R.
\]
Since $U_Q(t_Q)<R$ for each cube $Q\in\mathscr{D}$, adjoining cubes to~$\mathcal{Q}$ 
one by one preserves the hypothesis $\sum_{Q\in\mathcal{Q}} U_Q(t_Q)\le 2R$. 
Thus~\eqref{est:b-one-family} holds for every \emph{finite} family 
$\mathcal{Q}\in\mathfrak{S}$. Passing to arbitrary families in~$\mathfrak{S}$ 
by taking the supremum over finite subfamilies yields $C_b\le\lambda^\gamma C_2<\infty$.

We now prove the desired inequality~\eqref{ineq:Jensen}. Fix 
$Q\in\mathscr{D}$ and $g\in L^{p(\cdot)}(\Omega,\mu)$ with $\|g\|_{p(\cdot)}\le1$ 
and $|g|_Q>0$. Since $1<p_-\le p_+<\infty$, the function $t\mapsto I_Q(t)$ is continuous 
and strictly increasing from~$0$ to~$\infty$, so there exists a unique $t>0$ such that 
$I_Q(t)=|g|_Q$. By H\"older's inequality, the unit sphere property 
$\int_Q \|\mathbf{1}_Q\|_{p'(\cdot)}^{-p'(y)}\,d\mu(y)=1$, and our choice of~$\lambda$,
\[
I_Q(t)\le \frac{2}{\mu(Q)}\,\|g\|_{p(\cdot)} \|\mathbf{1}_Q\|_{p'(\cdot)}
\le 2\,I_Q \bigg(\frac{1}{\|\mathbf{1}_Q\|_{p'(\cdot)}}\bigg)
\le I_Q \bigg(\frac{\lambda}{\|\mathbf{1}_Q\|_{p'(\cdot)}}\bigg).
\]
The monotonicity of $I_Q$ implies
\begin{equation}\label{t:upper-est}
t\le\frac{\lambda}{\|\mathbf{1}_Q\|_{p'(\cdot)}}\le 
\lambda\,\max\bigg\{1,\frac{1}{\|\mathbf{1}_Q\|_{p'(\cdot)}^r}\bigg\}.
\end{equation}
Moreover, using Young's inequality 
$2t\,|g(x)|\le t^{p'(x)}+|2g(x)|^{p(x)}$, we deduce from the identity $I_Q(t)=|g|_Q$ 
that
\begin{align}\label{ineq:Young}
\int_Q t^{p'(y)}\,d\mu(y)&=
2t \int_Q |g(x)|\,d\mu(x) - \int_Q t^{p'(y)}\,d\mu(y)
\nonumber
\\
&\le 2^{p_+} \int_Q |g(x)|^{p(x)}\,d\mu(x).
\end{align}
We distinguish two cases.

\paragraph*{Case 1: $V_Q(t)\le\kappa\,U_Q(t)$.}
Raising both sides to the power~$r$ and applying the Young estimate~\eqref{ineq:Young} immediately yields
\[
\int_Q (|g|_Q)^{p(x)}\,d\mu(x)
=\int_Q I_Q(t)^{p(x)}\,d\mu(x)
\le \kappa^r\int_Q t^{p'(x)}\,d\mu(x)
\le C\int_Q |g(x)|^{p(x)}\,d\mu(x),
\]
where $C:=2^{p_+}\kappa^r>1$. This is inequality~\eqref{ineq:Jensen} without an error 
term.

\paragraph*{Case 2: $V_Q(t)>\kappa\,U_Q(t)$.}
Then $t<\min\big\{1,\|\mathbf{1}_Q\|_{p'(\cdot)}^{-r}\big\}$. 
Otherwise, the parameter 
\[
\tau:=\min\bigg\{t,\max\bigg\{1,\frac{1}{\|\mathbf{1}_Q\|_{p'(\cdot)}^r}\bigg\}\bigg\}
\]
lies in the characteristic window of Lemma~\ref{le:window}, which yields
$V_Q(\tau)\le C_1^\frac1r\,U_Q(\tau)$. Since we also have 
$\tau\le t\le \lambda\tau$ by~\eqref{t:upper-est}, this leads to the contradiction
\[
V_Q(t)\le \lambda^\gamma\,V_Q\Big(\frac{t}{\lambda}\Big) 
\le \lambda^\gamma\,V_Q(\tau)
\le \lambda^\gamma C_1^{\frac{1}{r}}\,U_Q(\tau)\le \kappa\,U_Q(t).
\]
Thus $t<1$, which corresponds to $|g|_Q=I_Q(t)<1$. Furthermore, 
$t^\frac1r<\|\mathbf{1}_Q\|_{p'(\cdot)}^{-1}$ gives
\[
\int_Q \big(t^\frac1r\big)^{p'(x)}\,d\mu(x)\le1.
\]
Applying the reverse H\"older property for~$p'(\cdot)$ from Lemma~\ref{le:RHI}(a) 
yields $U_Q(t)\le R$. Thus $t\in\mathcal{E}_Q$, which implies $t\le t_Q$ and 
$V_Q(t)\le V_Q(t_Q)=b(Q)$. Since $|g|_Q<1$, we have 
\[
\fint_Q (|g|_Q)^{p(x)}\,d\mu(x)\le(|g|_Q)^{p_-},
\]
and therefore with $\delta:=p_-(1-\frac{1}{r})>0$,
\[
\int_Q (|g|_Q)^{p(x)}\,d\mu(x)
=V_Q(t)\bigg(\fint_Q (|g|_Q)^{p(x)}\,d\mu(x)\bigg)^{1-\frac{1}{r}}
\le b(Q)\,(|g|_Q)^\delta\,\mathbf{1}_{(0,1)}(|g|_Q).
\]
Combining both cases completes the proof of~\eqref{ineq:Jensen}.
\end{proof}

\section{Proof of the Averaging Criterion}
\label{sec:proof}

Equipped with the self-improving property of the class 
$\mathcal{A}^\mathscr{D}(\Omega,\mu)$ and the Jensen-type inequality of 
Theorem~\ref{th:jensen}, we now prove Theorem~\ref{th:main}.
Our argument builds on a Calder\'on--Zygmund decomposition and 
a level-dependent truncation of the function.
\begin{theorem}
Let $(\Omega,d,\mu)$ be an \emph{unbounded} space of homogeneous type \emph{with a 
Borel-semiregular measure}~$\mu$, and let~$\mathscr{D}$ be a dyadic system on~$\Omega$. 
Suppose $p(\cdot)\in\mathcal{P}(\Omega,\mu)$ satisfies $1<p_-\le p_+<\infty$. 
Then~$M$ is bounded on~$L^{p(\cdot)}(\Omega,\mu)$ 
\emph{if and only if} $p(\cdot)\in\mathcal{A}^\mathscr{D}(\Omega,\mu)$.
\end{theorem}
\begin{proof}
Necessity is immediate from the pointwise inequality~\eqref{ineq:necessity}.
For sufficiency, assume $p(\cdot)\in\mathcal{A}^\mathscr{D}(\Omega,\mu)$. 
By Lemma~\ref{le:M-MD}(a), it suffices to fix $t\in\{1,\ldots,K\}$ and 
find a constant $c\ge1$ such that 
$\|M^{\mathscr{D}^t}f\|_{p(\cdot)}\le c$ for all 
$f\in L^{p(\cdot)}(\Omega,\mu)$ with $\|f\|_{p(\cdot)}\le1$. Applying the Fatou 
property to the truncations $f_n:=\min\{|f|,n\}\,\mathbf{1}_{B_n}$ over an expanding 
sequence of balls~$B_n$, we may assume without loss of generality that $f\ge0$ is 
bounded and supported in a ball, so that $f\in L^1(\Omega,\mu)$.

Lemma~\ref{le:RHI}(b) yields a scaling factor $s=\frac1r\in(\frac{1}{p_-},1)$ 
such that $q(\cdot):=sp(\cdot)$ satisfies $1<q_-\le q_+<\infty$ and belongs 
to~$\mathcal{A}^\mathscr{D}(\Omega,\mu)$. By the Calder\'on--Zygmund
Lemma~\ref{le:CZ}, for each $k\in\mathbb{Z}$, the level set
\[
\Omega_k:=\left\{x\in\Omega\ :\ M^{\mathscr{D}^t}f(x)>2\cdot 2^{ks}\right\}
=\bigcup_{Q\in\mathcal{Q}_k} Q
\]
decomposes into a family $\mathcal{Q}_k\in\mathfrak{S}$ of pairwise disjoint 
cubes $Q\in\mathscr{D}^t$ satisfying $|f|_Q>2\cdot 2^{ks}$. Together 
with the pointwise bound
$M^{\mathscr{D}^t}f\le 2\cdot 2^{(k+1)s}$ on $\Omega_k \setminus\Omega_{k+1}$,
this implies
\begin{align}\label{ineq:modular-rep}
\varrho_{p(\cdot)}(M^{\mathscr{D}^t}f)
&=\sum_{k\in\mathbb{Z}}\int_{\Omega_k\setminus\Omega_{k+1}}
M^{\mathscr{D}^t}f(x)^{p(x)}\,d\mu(x)
\nonumber
\\
&\le 2^{(1+s)\,p_+}
\sum_{k\in\mathbb{Z}} \sum_{Q\in\mathcal{Q}_k} \int_Q 2^{kq(x)}\,d\mu(x).
\end{align}

Define $f_k:=f\,\mathbf{1}_{\{f>2^{ks}\}}$. The three functions $f$, $f_k$,
and $f-f_k=f\,\mathbf{1}_{\{f\le2^{ks}\}}$ are nonnegative. For each cube 
$Q\in\mathcal{Q}_k$, we have $|f|_Q>2\cdot 2^{ks}$ and $|f-f_k|_Q\le 2^{ks}$, and thus
\[
|f_k|_Q=|f|_Q-|f-f_k|_Q>2^{ks}.
\]
We can therefore define the rescaled functions 
\[
g_{k,Q}:=\frac{2^k}{|f_k|_Q}\,f_k
\]
with average $|g_{k,Q}|_Q=2^k$. Since $\frac{2^k}{|f_k|_Q}<2^{k(1-s)}$, for all 
$x\in\Omega$ with $f(x)>0$ we have
\begin{align}\label{est:integrand}
g_{k,Q}(x)^{q(x)}&\le\big[2^{k(1-s)}\,f(x)\big]^{q(x)}\,\mathbf{1}_{\{f>2^{ks}\}}(x)
\nonumber
\\
&=f(x)^{p(x)}\,\Big(\frac{2^{ks}}{f(x)}\Big)^{(1-s)\,p(x)}\,
\mathbf{1}_{\{f>2^{ks}\}}(x)
\nonumber
\\
&\le f(x)^{p(x)}\,\Big(\frac{2^{ks}}{f(x)}\Big)^{(1-s)\,p_-}\,
\mathbf{1}_{\{f>2^{ks}\}}(x).
\end{align}
In particular, $g_{k,Q}(x)^{q(x)}\le f(x)^{p(x)}$, which implies 
$\varrho_{q(\cdot)}(g_{k,Q})\le\varrho_{p(\cdot)}(f)\le1$ and hence 
$\|g_{k,Q}\|_{q(\cdot)}\le1$. Applying the Jensen-type inequality from 
Theorem~\ref{th:jensen} with exponent~$q(\cdot)$ to~$g_{k,Q}$ yields constants 
$C>1$ and $\delta>0$ and an error function~$b$ with $C_b<\infty$---all independent 
of~$f$, $k$, $Q$ and~$t$---such that
\[
\int_Q 2^{kq(x)}\,d\mu(x) 
\le C\int_Q g_{k,Q}(x)^{q(x)}\,d\mu(x)
+b(Q)\,2^{k\delta}\,\mathbf{1}_{(0,1)}(2^k).
\]

Summing over the disjoint cubes $Q\in\mathcal{Q}_k$ and over $k\in\mathbb{Z}$, and 
substituting the pointwise bound~\eqref{est:integrand}, we obtain
\begin{align}\label{est:final-sum}
\sum_{k\in\mathbb{Z}}\sum_{Q\in\mathcal{Q}_k} &\int_Q 2^{kq(x)}\,d\mu(x) 
\nonumber
\\
&\le C\sum_{k\in\mathbb{Z}}\sum_{Q\in\mathcal{Q}_k} \int_Q g_{k,Q}(x)^{q(x)}\,d\mu(x)
+\sum_{k=-\infty}^{-1} 2^{k\delta}\sum_{Q\in\mathcal{Q}_k}b(Q)
\nonumber
\\
&\le C\int_{\{f>0\}} f(x)^{p(x)}
\bigg[
\sum_{k\in\mathbb{Z}} 
\Big(\frac{2^{ks}}{f(x)}\Big)^{(1-s)\,p_-}\,
\mathbf{1}_{\{f>2^{ks}\}}(x)
\bigg]\,d\mu(x)+\frac{C_b}{2^\delta-1}
\nonumber
\\
&\le\frac{C}{1-2^{-s(1-s)\,p_-}}+\frac{C_b}{2^\delta-1}<\infty.
\end{align}
The last inequality holds by evaluating 
the bracketed sum for each fixed $x\in\Omega$ 
with $f(x)>0$. Letting $k_0\in\mathbb{Z}$ be the largest integer such that 
$f(x)>2^{k_0 s}$, the sum in brackets reduces to
\[
\sum_{k=-\infty}^{k_0} \Big(\frac{2^{ks}}{f(x)}\Big)^{(1-s)\,p_-}
\le \sum_{k=-\infty}^{k_0} \big(2^{-s(1-s)\,p_-}\big)^{k_0-k}
=\frac{1}{1-2^{-s(1-s)\,p_-}},
\]
which we use together with $\varrho_{p(\cdot)}(f)\le1$. 
Combining~\eqref{ineq:modular-rep} with~\eqref{est:final-sum} gives 
\[
\varrho_{p(\cdot)}(M^{\mathscr{D}^t}f)
\le 2^{(1+s)\,p_+}\bigg(\frac{C}{1-2^{-s(1-s)\,p_-}}+\frac{C_b}{2^\delta-1}\bigg)=:c,
\]
where $c\ge1$ is independent of~$f$ and~$t$. By 
Theorem~\ref{th:toolbox}(e), we conclude $\|M^{\mathscr{D}^t}f\|_{p(\cdot)}\le c$.
\end{proof}

\bibliographystyle{abbrv}
\bibliography{AS-bib}

\end{document}